\documentclass[12pt]{amsart}
\usepackage{amsmath, amsthm, amssymb, color}
\usepackage[colorlinks=true,linkcolor=blue,urlcolor=blue]{hyperref}
\usepackage[top=25mm, bottom=25mm, left=25mm, right = 25mm]{geometry}
\usepackage{graphicx}
\usepackage{mathtools}
\usepackage{enumerate}
\usepackage{verbatim}
\usepackage{tikz,tikz-cd,tikz-3dplot}
\usepackage{amssymb}
\usepackage{makecell}
\usepackage{array}
\usepackage{enumitem}
\usepackage{multirow}
\usepackage{afterpage}
\usepackage{algorithm}
\usepackage{algpseudocode}
\usepackage{booktabs}
\usepackage{caption}
\usepackage[labelformat=simple]{subcaption}

\allowdisplaybreaks

\newtheorem{theorem}{Theorem}[section]
\newtheorem{proposition}[theorem]{Proposition}
\newtheorem{lemma}[theorem]{Lemma}
\newtheorem{corollary}[theorem]{Corollary}
\newtheorem{conjecture}[theorem]{Conjecture}

\theoremstyle{definition}
\newtheorem{definition}[theorem]{Definition}
\newtheorem{problem}[theorem]{Problem}
\newtheorem{remark}[theorem]{Remark}

\newtheorem{example}[theorem]{Example}

\newcommand{\PP}{\mathbb{P}}
\newcommand{\RR}{\mathbb{R}}

\newcommand{\CC}{\mathbb{C}}

\newcommand{\M}{\mathcal{M}}
\newcommand{\PGL}{\mathrm{PGL}}
\newcommand{\sgn}{\mathrm{sgn}}

\definecolor{dgreen}{HTML}{026a10}
\definecolor{dviolet}{HTML}{9109E3}
\definecolor{dorange}{HTML}{e55700}

\title{Dihedral sign patterns in $\M_{0,n}$}
\author{Veronica Calvo Cortes}
\address{Veronica Calvo Cortes (MPI MiS)}
\email{veronica.calvo@mis.mpg.de}

\author{Hannah Tillmann-Morris}
\address{Hannah Tillmann-Morris (MPI MiS)}
\email{hannah.tillmann@mis.mpg.de}

\date{}

\begin{document}

\maketitle


\begin{abstract}
    The connected components of $\M_{0,n}(\RR)$ are in bijection with the $(n-1)!/2$ dihedral orderings of $[n]$. They are all isomorphic. 
    We construct monomial maps between them, and use these maps to prove a conjecture of Arkani-Hamed, He, and Lam in the case of $\M_{0,n}$. Namely, we provide a bijection between connected components and sign patterns that are consistent with the extended $u$-relations for the dihedral embedding.
\end{abstract}

\section{Introduction}
The configuration space $\M_{0,n}$ of $n$ distinct points on the projective line $\PP^1$ up to projective equivalence is a well-studied object. It is combinatorially and geometrically rich, and together with its Deligne-Mumford-Knudsen compactification $\overline{\M}_{0,n}$, it has proven to be a key example in the emerging field of positive geometry \cite{LamNotes, whatis, scatteringcurves}.  We focus mainly on the combinatorics of the real points $\M_{0,n}(\RR)$. This variety has $(n-1)!/2$ connected components, each corresponding to an ordering of the $n$ points on $\PP^1(\RR) \cong S^1$. Our main result, Theorem~\ref{thm:mainResult}, is an alternative characterization of the connected components using sign patterns.

\medskip

To consider sign patterns, we need explicit coordinates for $\M_{0,n}(\RR)$. We will work with two different descriptions of this space: first as a hyperplane arrangement complement and then via the cross-ratio parametrizations. We can represent $z \in \M_{0,n}$ as a $2 \times n$ matrix whose columns correspond to the marked points
\begin{equation}\label{eq:matrix_rep}
    M_z = \begin{bmatrix}
        1 & 1 &\cdots & 1\\
        z_1 & z_2 &\cdots &z_n
    \end{bmatrix}.
\end{equation}
Acting by an element of $\PGL(2)$, we place $z_1=0, z_2=1$ and $z_n=\infty$ (i.e. $[1: z_n]=[0:1]$). 

\subsection{Hyperplane arrangement} \label{subsec: hypArr}
We can view $\M_{0,n}$ as the complement of a hyperplane arrangement in $\CC^{n-3}$. As we have fixed three of our $n$ marked points, we only need $n-3$ coordinates $z_3,\ldots, z_{n-1}$. Each $z_i$ must be distinct from $0,1,\infty$ as well as being pairwise distinct. Thus the arrangement we consider consists of the hyperplanes
\begin{align}\label{eq:M0nhyperplanes}
    H_{0,i} &= \{z \in \CC^{n-3}: z_i=0\} && \text{ for } i=3,\ldots, n-1 \nonumber\\
    H_{1,i} &= \{z \in \CC^{n-3}: z_i=1\} && \text{ for } i=3,\ldots, n-1\\
    H_{i,j} &= \{z \in \CC^{n-3}: z_i=z_j\} && \text{ for } 3\leq i < j \leq  n-1 \nonumber.
\end{align}
Note that to study the real points of $\M_{0,n}$, we may view it as the complement of the hyperplanes in \eqref{eq:M0nhyperplanes} in $\RR^{n-3}$. The regions of the real hyperplane arrangement $\mathcal{A}$ encode the signs of the differences $z_i-z_j$ for $1\leq i < j \leq n-1$ (see more in \cite[Section 19]{LamMatroidsAmp}).

\subsection{Dihedral coordinates}
We are interested in understanding sign patterns in a specific parametrization of $\M_{0,n}$. Recall that for any distinct four indices $i,j,k,l \in [n]$ the cross-ratio of their corresponding points in $\PP^1$ is given by
\[
    [ij|kl]:= \frac{(z_i-z_k)(z_j-z_l)}{(z_i-z_l)(z_j-z_k)}=\frac{p_{ik}(M_z)p_{jl}(M_z)}{p_{il}(M_z)p_{jk}(M_z)},
\]
where the $p$'s denote the Plücker coordinates of $M_z$. These define an embedding 
\[\iota: \M_{0,n} \to (\PP^1\setminus\{0,1,\infty\})^{\binom{n}{4}} \qquad z \mapsto \left([ij|kl]: {i,j,k,l} \subset [n]\right).\]
The image of $\iota$ is cut out by the cross-ratio relations \eqref{eq: crossratiorels}:
\begin{align}\label{eq: crossratiorels}
    &[ij|kl]= [ji|lk]=[kl|ij]=[lk|ji]  &&[ij|lk]= [ji|kl] = \frac{1}{[ij|kl]} \\
    &[ik|jl]=1-[ij|kl] && [ij|kl]=[ij|km][ij|ml]  \qquad \forall i,j,k,l,m \in [n]. \nonumber
\end{align}

The closure of $\M_{0,n}$ in $(\PP^1)^{\binom{n}{4}}$ under the embedding $\iota$ is a smooth projective variety $\overline{\M}_{0,n}$. This is the well-studied DMK compactification~\cite{Knudsen}, the moduli space of stable rational curves. The boundary of $\overline{\M}_{0,n}$ satisfies a factorization property: each of  $2^{n-1}-n-1$ connected components of the boundary is isomorphic to a product of smaller moduli spaces $\overline{\M}_{0,n_1} \times \overline{\M}_{0,n_2}$, where $n_1 + n_2 = n + 2$. 
This stratification property of $\overline{\M}_{0,n}$ leads to the factorization of the integrals appearing in the study of positive geometries and scattering amplitudes \cite{stringy}.\\

In this paper we study the embedding of $\M_{0,n}$ into an algebraic torus of lower dimension, by considering a subset of the cross-ratios motivated by the study of its positive part.

\begin{definition}\label{def: dihedral coords}
    The \emph{dihedral coordinates} of $\M_{0,n}$ are the following cross-ratios:
\[
    u_{ij} = [i, i+1 | j+1, j], \qquad \text{ for } ij  \text{ a chord of the $n$-gon.} 
\]
\end{definition}

We can consider different sets of dihedral coordinates depending on distinct labelings of the vertices of the $n$-gon. These labelings correspond to elements in $S_n/D_{2n}$, known as \emph{dihedral orderings} of $[n]$. Note that the above definition of the dihedral coordinate $u_{ij}$ depends on the vertices of the $n$-gon being labeled with the standard dihedral ordering $1,2,3, \ldots, n$ -- the numbers $i+1$ and $j+1$ are the labels of the vertices to the right of $i$ and $j$ respectively. The positive part $\M_{0,n}(\RR_{>0})$ is defined to be the connected component of $\M_{0,n}(\RR)$ where the $n$ points are ordered on $\PP^1(\RR) \cong S^1$ with respect to this ordering. This is exactly the component of $\M_{0,n}(\RR)$ where these dihedral coordinates are positive.

\subsection{U-relations}\label{subsec: 1.3}
The dihedral coordinates satisfy nice combinatorial relations which cut out $\M_{0,n}$ set-theoretically. For $ij$ a chord of the $n$-gon we consider the (primitive) $u$-relation
\[
    R_{ij} := u_{ij} + \prod_{kl \nsim ij} u_{kl}-1
\]
where $kl \nsim ij$ means that the chord $kl$ crosses the chord $ij$. 

\begin{proposition}\cite[Section 2.2]{Brown09}
    The map $\M_{0,n} \to (\CC^*)^{\binom{n}{2}-n}$ given by the dihedral coordinates is a closed embedding. Moreover, the ideal of the image is $\langle R_{ij}: ij$ chord of $n$-gon$\rangle$.
\end{proposition}
\begin{proof}
    Let $I=\langle R_{ij}: ij$ chord of $n$-gon$\rangle \subset \CC[u_{ij}^\pm]$ be the ideal generated by the $u$-relations (for the standard dihedral ordering). We check that the vanishing of $I$ is the image of the map. One containment follows from a direct computation: replacing the $u_{ij}$ with the corresponding cross-ratio. For the other containment we take a point in $\mathbb{V}(I)$ and build the corresponding configuration of $n$ points in $\PP^1$. Recall that without loss of generality we are assuming $z_1=0, z_2=1$ and $z_n=\infty$. Then for $i\in \{3,\ldots,n-1\}$ we have 
    \[
        u_{in}=[i \ i+1 | 1\, n] = \frac{(z_i-z_1)\det \begin{bmatrix}
            1 & 0 \\ z_i & 1
        \end{bmatrix}}{\det \begin{bmatrix}
            1 & 0 \\ z_{i+1} & 1
        \end{bmatrix}(z_{i+1}-z_1)} = \frac{z_i}{z_{i+1}}.
    \]
    Then, we can set $z_3=u_{2n}^{-1}$, $z_4=(u_{3n}u_{2n})^{-1}$ and so on until $z_{n-1}=(u_{n-2\,n}\cdots u_{2n})^{-1}$. The fact that the $u$-relations vanish guarantees that this yields $n$ distinct points. 
\end{proof}

Among the elements of the ideal generated by the $u$-relations (both in the Laurent and usual polynomial ring) there is another special set of relations (see \cite[Proposition 1.8]{LamNotes}): for any partition of $[n]$ into four cyclic intervals $A=(a,a+1,\ldots,b-1), B=(b,b+1,\ldots,c-1), C=(c,c+1,\ldots,d-1)$ and $D=(d,d+1,\ldots,a-1)$ we consider the \emph{extended $u$-relation}
\begin{equation}\label{eq: ext u eq}
    R_{A,B,C,D} := \prod_{i \in A,\; j \in C} u_{ij}+\prod_{k \in B, \; l \in D} u_{kl}-1.
\end{equation}
We also denote the extended $u$-relation $R_{A,B,C,D}$ by 
$
R_{a, \ldots, b-1 | b, \ldots, c-1 | c, \ldots, d-1| d, \ldots, a-1}
$
if convenient.

\begin{definition}\label{def: sign patterns}
    Let $\RR_s$ denote the orthant in $\RR^{\binom{n}{2}-n}$ determined by $s \in \{-,+\}^{\binom{n}{2}-n}$. We call $s$ a \emph{consistent sign pattern} if it does not contradict the extended $u$-relations \eqref{eq: ext u eq}. We call it a \emph{realizable sign pattern} if 
    $\M_{0,n}(\RR_s):=\RR_s \cap \M_{0,n}(\RR)$ is non-empty.
\end{definition}

We are now ready to state our main result. This is \cite[Conjecture 11.1]{ClusterConfigSpaces} for type $A$.

\begin{theorem}\label{thm:mainResult}
    $\M_{0,n}(\RR_s)$ is non-empty if and only if $s$ is consistent. Equivalently, all consistent sign patterns are realizable. Moreover, if $\M_{0,n}(\RR_s) \neq \emptyset$ then it is connected and isomorphic as a semi-algebraic set to $\M_{0,n}(\RR_{>0})$.
\end{theorem}

As suggested by Arkani-Hamed, He, and Lam \cite{ClusterConfigSpaces}, consistent sign patterns are analogous to uniform oriented matroids. We hope that alternative combinatorial characterizations of consistency (as for matroids) will yield a more elegant proof of Theorem~\ref{thm:mainResult}. Our approach directly uses Definition~\ref{def: sign patterns}, and hence is long and technical.

\begin{example}[$n=6$]
    There are nine dihedral coordinates $(u_{13},u_{14},u_{15},u_{24}, u_{25},u_{26},u_{35}$, $u_{36},u_{46})$ on $\M_{0,6}$. The nine primitive $u$-relations are
    \begin{align*}
        u_{13} + u_{24}u_{25}u_{26}-1, && u_{24}+u_{35}u_{36}u_{13}-1, && u_{35}+u_{46}u_{14}u_{24}-1,\\
    u_{46}+u_{15}u_{25}u_{35}-1, && u_{15}+u_{26}u_{36}u_{46}-1, && u_{26}+u_{13}u_{14}u_{15}-1,\\
 u_{14}+u_{25}u_{26}u_{35}u_{36}-1, && u_{25}+u_{36}u_{13}u_{46}u_{14}-1, &&u_{36}+u_{14}u_{24}u_{15}u_{25}-1,
    \end{align*}
    and there are six further extended $u$-relations:
    \begin{align}
        u_{14}u_{24} + u_{35}u_{36}-1, && u_{25}u_{35} + u_{46}u_{14} - 1, && u_{36}u_{46} + u_{15}u_{25} -1, \nonumber\\
 u_{14}u_{15} + u_{26}u_{36} - 1, && u_{25}u_{26} + u_{13}u_{14} - 1, && u_{36}u_{13} + u_{24}u_{25}-1. \label{eq: n=6extended}
    \end{align}
    The following $14$ sign patterns are consistent with the nine primitive $u$-relations, but are not consistent in the sense of Definition~\ref{def: sign patterns}, as they contradict the six extended $u$-relations~\eqref{eq: n=6extended}.
    \begin{align*}
    (-,-,+,-,+,-,-,+,+) && (-,-,+,+,-,-,+,+,+) &&(-,+,-,-,+,-,-,+,-)\\
  (-,+,-,-,+,-,+,-,+) &&(-,+,-,+,-,-,+,+,-) &&(-,+,+,-,-,+,-,+,-)\\
  (-,+,+,-,-,+,+,-,+) &&(+,-,-,+,+,-,-,+,-) &&(+,-,-,+,+,-,+,-,+)\\
  (+,-,+,-,+,+,-,-,+) &&(+,-,+,+,-,+,-,+,-) &&(+,-,+,+,-,+,+,-,+)\\
  (+,+,-,-,+,+,-,-,-) &&(+,+,-,+,-,+,+,-,-)&&
    \end{align*}
    By Theorem~\ref{thm:mainResult}, there are precisely 60 consistent sign patterns.
    The nine primitive $u$-relations generate a prime ideal in the polynomial ring, as expected (see~\cite[Theorem 3.3]{ClusterConfigSpaces}). For $n = 5$, the five $u$-relations and 12 consistent sign patterns are listed in Example~\ref{ex: n=5}.
    \begin{table}[h!]
    \centering
    \begin{tabular}{|c|c|c|c|c|}
        \hline
         $n$ & 5 & 6 & 7 & 8 \\
         \hline
          primitive & 12 & 74 & 697 & 10\,180\\
         extended & 12 & 60 & 360 & 2520\\
        \hline
    \end{tabular}
    \caption{The number of sign patterns consistent with just the primitive $u$-relations, versus those consistent with all extended $u$-relations.}
    \label{tab: no consistent sign patterns for n}
\end{table}
\end{example}

\subsection{Notation}\label{subsec: notation} 
In what follows we will be constantly switching between different coordinate systems. We write $\M_{0,n}^\alpha$ when using the $u$-variables corresponding to the dihedral ordering $\alpha$; we call this a \emph{dihedral chart}. We abuse notation by also regarding $\alpha$ as a representative in $S_n$ of the coset in $S_n/D_{2n}$. We drop the superscript when $\alpha=\text{id} \in S_n$ is the standard dihedral ordering $1,2,\ldots, n-1, n$. We use the notation  
    \begin{equation}\label{eq: notation uij}
            u_{ij}^\alpha := [\alpha(i) \, \alpha(i+1) \, | \, \alpha(j+1) \, \alpha(j)]
    \end{equation}
for the cross-ratio corresponding to the chord between the vertices $\alpha(i)$ and $\alpha(j)$ on the $n$-gon labeled by $\alpha$. These are the variables in the coordinate ring $\CC[\M_{0,n}^\alpha]$.

A \emph{sign pattern} $s$ is the choice of an orthant in $\RR^{\binom{n}{2} - n}$ and we denote it $\RR_s$. We denote the intersection $\M_{0,n}^\alpha(\RR) \cap \RR_s$ by $\M_{0,n}^\alpha(\RR_s)$. When $s =(+,\ldots,+)$ we replace $\RR_s$ by $\RR_{>0}$.

\subsection{Outline} In Section~\ref{sec:2} we introduce the main ingredient of the paper: monomial maps which provide isomorphisms between the different dihedral charts of $\M_{0,n}$. This allows us to understand the connected components of $\M_{0,n}(\RR)$ in terms of realizable sign patterns. In Section~\ref{sec:3} we prove Theorem~\ref{thm:mainResult} by providing an algorithm that associates a dihedral ordering to a consistent sign pattern. 
We implemented this algorithm; the code is available at our MathRepo website~\cite{M2Code}.
This section is the technical heart of the paper and is subdivided into four parts. 
In Section~\ref{sec:4} we discuss two interesting further directions: the cluster configuration space of type $C_n$ and the commutative algebra of $u$-relations. For examples and some proofs we used the computer algebra system Macaulay2 \cite{M2}. The scripts are also available at~\cite{M2Code}.

\section{Realizable sign patterns}\label{sec:2}

Our goal in this section is to prove that each realizable sign pattern corresponds to a unique connected component of $\M_{0,n}(\RR)$.
The key step is to build the monomial transformations which allow us to change between the different dihedral charts of $\M_{0,n}$. We first introduce the analogy to simple transpositions for dihedral orderings.

\begin{definition}
    Given a dihedral ordering represented by $\sigma \in S_n$, we define its \textit{adjacent transpositions} as the transpositions in $S_n$ permuting two vertices which are adjacent in the dihedral ordering. Explicitly, we denote the adjacent transposition which swaps the vertices $\sigma(k)$ and $\sigma(k+1)$ by $\tau_k^\sigma$. For ease of notation, we also refer to it with its cycle notation -- i.e. $(\sigma(k) \; \sigma(k+1)) \in S_n$ -- when the dihedral ordering $\sigma$ is clear from the context.
\end{definition}

The cross-ratio relations \eqref{eq: crossratiorels} allow us to define special monomial maps between dihedral charts of $\M_{0,n}$ whose dihedral orderings differ only by an adjacent transposition. 

\begin{definition}\label{def: monomial map for general adjacent transposition}
Let $\sigma \in S_n$ and take $\alpha = \tau_k^\sigma \sigma$ for some $k \in [n]$. We define the monomial map $\phi_k^\sigma: \CC[\M_{0,n}^\alpha] \rightarrow \CC[\M_{0,n}^\sigma]$ in generators as
    \begin{align}\label{eq: phi alpha k}
        u_{ij}^\alpha & \longmapsto u_{ij}^\sigma \quad & \text{if } i,j \notin \{k-1, k, k+1\}\nonumber\\
        u_{i \, k-1}^\alpha & \longmapsto u_{i \, k-1}^\sigma u_{i \, k}^\sigma \quad & \text{if } i \notin \{k-1, k, k+1\}\nonumber\\
        u_{i \, k+1}^\alpha & \longmapsto u_{ik}^\sigma u_{i\,k+1}^\sigma \quad & \text{if } i \notin \{k-1, k, k+1\}\\
        u_{ik}^\alpha & \longmapsto (u_{ik}^\sigma)^{-1}  \quad & \text{if } i \notin \{k-1, k, k+1\}\nonumber\\
        u_{k-1 \, k+1}^\alpha  & \longmapsto - \, u_{k-1\,k+1}^\sigma \prod_{i \notin \{k, k\pm 1\}} (u_{ik}^\sigma)^{-1}. \nonumber
    \end{align}
This defines a transition map between dihedral charts $\varphi_\alpha^\sigma: \M_{0,n}^\sigma \rightarrow \M_{0,n}^\alpha$.
\end{definition}

\begin{example}\label{ex: monomial map for (12)}
We illustrate how the map in Definition~\ref{def: monomial map for general adjacent transposition} comes from the cross-ratio relations~\eqref{eq: crossratiorels}. Take $\sigma=\text{id}$ and $k=1$, so that $\alpha = (12)$. We relabel the coordinates on $\M_{0,n}^{(12)}$ by $v_{\alpha(i) \, \alpha(j)} := u_{ij}^{\alpha}$. Then,
    \begin{align*}
        v_{ij} = u_{ij}^{(12)} & = [i \, i+1\,| \, j+1 \, j] = u_{ij} \quad & \text{for all } i,j \notin \{n,1,2\}\\
        v_{ni} = u_{ni}^{(12)} & = [n \, 2 | \, i+1 \, i] &\\
        & = [n \, 1 \, | \, i+1 \, i] [1 \, 2 \, | \, i+1 \, i] = u_{ni}u_{1i} \quad & \text{for all } i \notin \{n,1,2\}\nonumber\\
        v_{1i} = u_{2i}^{(12)} & = [1 \, 3 | \, i+1 \, i] &\\
        & = [1 \, 2 \, | \, i+1 \, i] [2 \, 3 \, | \, i+1 \, i] = u_{1i}u_{2i} \quad & \text{for all } i \notin \{n,1,2\}\\
        v_{2i} = u_{1i}^{(12)} & = [2 \, 1 \, | \, i+1 \, i] & \\
        & = [1 \, 2 \, | \, i+1 \, i]^{-1} = u_{1i}^{-1} \quad & \text{for all } i \notin \{n,1,2\}\\
        v_{n1} = u_{n2}^{(12)} & = [n \, 2 \, | \, 3 \, 1] & \\
        & = [1 \, 3 \, | \, n \, 2]^{-1} & \\
        & = \left( 1 - [1 \, n \, | \, 3 \, 2] \right)^{-1} & \\
        & = \frac{[n \, 1 \, | \, 3 \, 2]}{[n \, 1 \, | \, 3 \, 2] - 1} = \frac{u_{n2}}{u_{n2} - 1} = - \, u_{n2} \prod_{i \notin \{n,1,2\}} u_{1i}^{-1}.
    \end{align*}
\end{example}

\begin{lemma}\label{lem:monomial_map}
    Let $\alpha$ be a dihedral ordering of $[n]$. There is an invertible monomial transformation $\phi_\alpha:\CC[\M_{0,n}^\alpha] \to \CC[\M_{0,n}] $ providing a change of coordinates between dihedral charts.
\end{lemma}

\begin{proof}
    For any dihedral ordering $\alpha$ we can ``sort'' it to the standard ordering $1,2,\ldots,n$ by performing only adjacent transpositions. Composing the corresponding monomial maps from Definition~\ref{def: monomial map for general adjacent transposition} we obtain a map $\phi_\alpha: \CC[\M_{0,n}^{\alpha}] \to \CC[\M_{0,n}].$ Explicitly, if the sequence of adjacent transpositions that sorts $\alpha$ to the standard ordering is $\tau_{i_k}^{\alpha_k}\cdots \tau_{i_1}^{\alpha_1}$, where $\alpha_1=\alpha$ and $\alpha_{j+1}=\tau_{i_j}^{\alpha_j}\cdots \tau_{i_1}^{\alpha_1}\alpha$ for all $1 \leq j\leq k-1$, then 
    $
    \phi_\alpha:=\phi_{i_k} \circ \phi_{i_{k-1}}^{\alpha_k}\circ \cdots \circ \phi_{i_1}^{\alpha_2}.
    $
\end{proof}

\begin{remark}\label{rem: varphi alpha}
   We let $\varphi_\alpha: \M_{0,n} \rightarrow \M_{0,n}^\alpha$ denote the isomorphism between the two different embeddings of $\M_{0.n}$ which is given by the monomial map $\phi_\alpha$ (defined in Lemma~\ref{lem:monomial_map}). Viewing $\alpha$ as a permutation in $S_n$, the map $\varphi_\alpha$ is an automorphism of $\M_{0,n}$ encoding the action of permuting the marked points.  
\end{remark}

\begin{example}[$n=5$]\label{ex: n=5}
    We consider the surface $\M_{0,5}$ with the standard dihedral coordinates $u_{13}, u_{14}, u_{24},$ $u_{25}, u_{35}$. It has $12$ connected components corresponding to the dihedral orderings of $\{1,2,3,4,5\}$. Each of them corresponds to a fixed sign pattern in the $u$-variables, which we calculate by computing the image of each generator under $\phi_\alpha$.
    
    Explicitly, let $\alpha$ be the dihedral ordering $1\,4\,2\,5\,3$, and denote the $u$-coordinates for $\M_{0,5}^\alpha$ by $v_{12},v_{15},v_{23},v_{34},v_{45}$, following the convention set in Example~\ref{ex: monomial map for (12)}. We can sort $\alpha$ to the standard ordering via the sequence of adjacent transpositions $(23)(24)(13)$. This yields the map $\phi_\alpha = \phi_2\circ \phi_2^{32451} \circ \phi_1^{34251} : \CC[\M_{0,5}^{14253}] \to \CC[\M_{0,5}]$, given by
    \begin{align*}
        &v_{12} \mapsto -\frac{u_{14}}{u_{25}u_{35}} &&v_{15}\mapsto -\frac{u_{35}}{u_{14}u_{24}} &&v_{23}\mapsto -\frac{u_{25}}{u_{13}u_{14}} \\
        &v_{34}\mapsto -\frac{u_{13}}{u_{24}u_{25}} && v_{45}\mapsto -\frac{u_{24}}{u_{13}u_{35}}.
    \end{align*}
    Positive values on all the $v$'s can only be achieved with negative values on all the $u$'s. Hence, $\alpha$ corresponds to the sign pattern with all negative coordinates. Similar computations give a correspondence between dihedral orderings and sign patterns $\sgn(u_{13},u_{14},u_{24},u_{25},u_{35})$.
    \begin{align*}
            &1\,2\,3\,4\,5 \quad  (+,+,+,+,+) &&& 1\,5\,3\,2\,4 \quad  (-,-,+,+,+)\\
            &1\,3\,2\,4\,5 \quad  (-,+,+,+,+) &&& 1\,5\,2\,4\,3 \quad  (+,-,-,+,+)\\
            &1\,5\,2\,3\,4 \quad  (+,-,+,+,+) &&& 1\,4\,3\,5\,2 \quad  (+,+,-,-,+)\\
            &1\,2\,4\,3\,5 \quad  (+,+,-,+,+) &&& 1\,3\,5\,4\,2 \quad  (+,+,+,-,-)\\
            &1\,3\,4\,5\,2 \quad  (+,+,+,-,+) &&& 1\,3\,2\,5\,4 \quad  (-,+,+,+,-)\\
            &1\,2\,3\,5\,4 \quad  (+,+,+,+,-) &&& 1\,4\,2\,5\,3 \quad  (-,-,-,-,-)
    \end{align*}
    See \cite[Eq. 8]{delpezzo} for another perspective on these $12$ sign patterns.
\end{example}

Having set up the monomial transformations, we can now relate dihedral orderings to realizable sign patterns. The following basic combinatorial lemma is required for Proposition~\ref{prop: section 2}. The linear description of $\M_{0,n}$ as a hyperplane arrangement complement will be key.

\begin{lemma}\label{prop: realizable sign patterns}
    There are $(n-1)!/2$ realizable sign patterns.
\end{lemma}
\begin{proof}
We view $\M_{0,n}(\RR)=\RR^{n-3}\setminus \mathcal{A}$ as a hyperplane arrangement complement as in Section~\ref{subsec: hypArr}. 
    By \cite[p.63]{LamMatroidsAmp}, the number of regions of $\mathcal{A}$ is $(n-1)!/2$, so it suffices to show that the realizable sign patterns are in bijection with the regions of $\mathcal{A}$.
    Recall that the dihedral coordinates are defined by multiplying/dividing differences of the form $z_k-z_l$. 
    It is thus clear how to associate a realizable sign pattern to any region in $\mathcal{A}$.
    We now give the opposite direction. That is, we determine the collection of signs $P(i,j) := \text{sgn}(z_i - z_j)$ for $1 \leq i < j \leq n$, given a sign pattern $s$ in the $u_{ij}$'s.

    We proceed by double induction on $i$ and $j$.
    For $i=1, j=2$ we must have $P(1,2)=-$ as we have assigned $z_1=0$ and $z_2=1$. Now we fix $i=1$ and assume we have computed the sign of $P(1,j-1)$. Recall that 
    \[
        s(u_{j-1\,n}) = \sgn \frac{(z_{j-1}-z_1)(z_j-z_n)}{(z_{j-1}-z_n)(z_j-z_1)} = P(1,j-1)P(j,n)P(j-1,n)P(1,j).
    \]
    Note that for any $k=1,\ldots,n-1$ we have $P(k,n)=-$, as we have assigned $z_n = \infty$. Hence, $P(1,j)=P(1,j-1)s(u_{j-1\,n})$, concluding the base case for $i$. 
    Now, we let $i \geq 2$ and assume that we have computed $P(k, l)$ for all $k<i$ and $l>k$, as well as $P(i, j-1)$.
    Note that
    \[
        s(u_{i-1\,j-1}) = \sgn \frac{(z_{i-1}-z_{j})(z_{i}-z_{j-1})}{(z_{i-1}-z_{j-1})(z_i-z_j)} = P(i-1,j)P(i,j-1)P(i-1,j-1)P(i,j).
    \]
    By the induction hypothesis, this allows us to compute $P(i,j)$ from $s$. 
\end{proof}

\begin{proposition}\label{prop: section 2}
    There is a one-to-one correspondence between realizable sign patterns and dihedral orderings. This is given by $\varphi_\alpha(\M_{0,n}(\RR_{s}))=\M_{0,n}^\alpha(\RR_{>0})$, where $s$ is a realizable sign pattern and $\varphi_\alpha$ is the map from Remark~\ref{rem: varphi alpha} associated to a dihedral ordering $\alpha$. Moreover, $\M_{0,n}(\RR_s)$ is a connected component of $\M_{0,n}(\RR)$.
\end{proposition}

\begin{proof}
     We use the notation from Remark~\ref{rem: varphi alpha} and denote $\varphi_\alpha$ for the geometric map and $\phi_\alpha$ for the induced map between coordinate rings. Lam's argument in \cite[Proposition 1.14]{LamNotes} shows that the images $\varphi_\alpha^{-1}(\M_{0,n}^{\alpha}(\RR_{>0}))$ are the connected components of $\M_{0,n}(\RR)$. Now, note that $\phi_\alpha^{-1}$ is also a monomial map and hence the whole image $\varphi_\alpha^{-1}(\M_{0,n}^{\alpha}(\RR_{>0}))$ is contained in a single orthant with sign pattern $s = s(\alpha)$. 
     In other words, this is a surjection from dihedral orderings onto realizable sign patterns.
     There are exactly $(n-1)!/2$ dihedral orderings on $[n]$, and so by Lemma~\ref{prop: realizable sign patterns}, this is a bijection.
\end{proof}

Although Proposition~\ref{prop: section 2} is very similar to our main result, extending the proof from realizable to consistent sign patterns (recall Definition~\ref{def: sign patterns}) is a tough problem in real algebraic geometry. 
The challenge is to solve the equations $(R_{ij} = 0)$ with prescribed signs.

\section{Dihedral orderings from sign patterns}\label{sec:3}
In this section we prove our main result, Theorem~\ref{thm:mainResult}, by producing a dihedral ordering from a consistent sign pattern $s$. Concretely, Algorithm~\ref{alg: main} outputs a dihedral ordering $\alpha \in S_n$ such that $\varphi_\alpha(\M_{0,n}(\RR_{s}))=\M_{0,n}^\alpha(\RR_{>0})$. We prove it terminates when $s$ is consistent.
\begin{algorithm}
\caption{Dihedral ordering from sign pattern}\label{alg: main}
\begin{algorithmic}
\Require $s$ is a consistent sign pattern
\State $t := s$ and $\alpha := \text{id} \in S_n$ 
\While{$t \neq (+,\ldots, +)$} \Comment{$(+,\ldots, +)$ is the all positive sign pattern}
    \State $u_{ab}^\alpha \gets$ ``shortest negative chord in $t$''
    \State $\alpha \gets \alpha=(\alpha(a+1)\,\alpha(b))\circ \alpha$
    \State $t \gets $ sign pattern of $\varphi_\alpha(\M_{0,n}(\RR_{s}))$ \Comment{$\varphi_\alpha$ is the monomial map from Lemma~\ref{lem:monomial_map}}
\EndWhile
\end{algorithmic}
\end{algorithm}

\noindent The implementation is available in~\cite{M2Code}. The following example illustrates what the algorithm does with a sign pattern in $\M_{0,10}$.

\begin{example}
    Let $s$ be the consistent sign pattern with five negative chords (in \textcolor{red}{red}) shown in the top left of Figure~\ref{fig:ExAlg}. The shortest negative chord is $u_{68}$ so we apply the transposition $(7 8)$. The new shortest negative chord is between vertices $7$ and $10$ -- in our notation $u_{8 \, 10}^{(7 8)}$ -- so we apply $(9 \, 10)$. The algorithm outputs $\alpha = 1\, 2\, 3\, 8\, 4\, 5\, 6\, 9\, 10\, 7$ after six iterations. Equivalently, $\alpha$ is represented by the sequence of transpositions $ (54)(68)(79)(46)(9 \, 10)(78)\in S_n$.
     \begin{figure}[h!]
  \centering
  \begin{subfigure}[c]{0.33\textwidth}
  	\includegraphics[width=\textwidth]{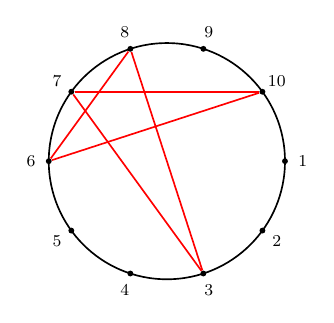}
  \end{subfigure}%
  \hfill %
  \begin{subfigure}[c]{0.33\textwidth}
  \centering
	\includegraphics[width=\textwidth]{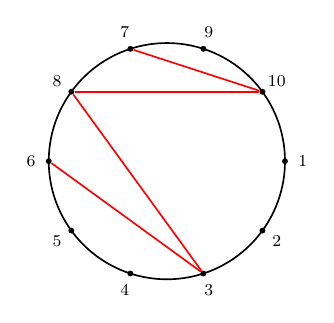}
  \end{subfigure}%
  \hfill%
  \begin{subfigure}[c]{0.33\textwidth}
  \centering
  	\includegraphics[width=\textwidth]{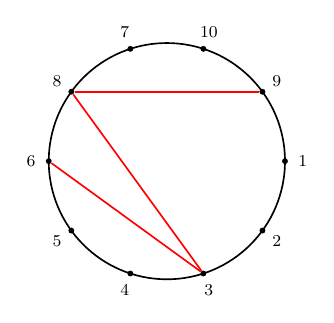}
  \end{subfigure}
  \newline
  \begin{subfigure}[c]{0.33\textwidth}
  	\includegraphics[width=\textwidth]{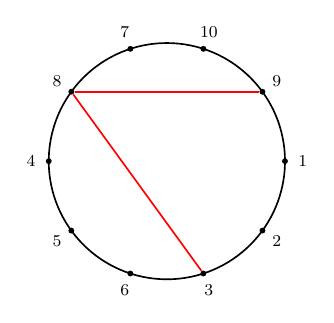}
  \end{subfigure}%
  \hfill %
  \begin{subfigure}[c]{0.33\textwidth}
  \centering
	\includegraphics[width=\textwidth]{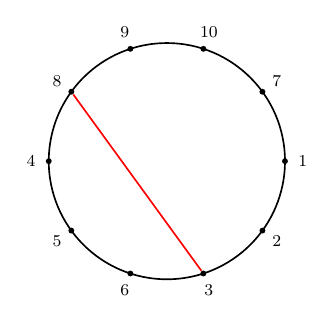}
  \end{subfigure}%
  \hfill%
  \begin{subfigure}[c]{0.33\textwidth}
  \centering
  	\includegraphics[width=\textwidth]{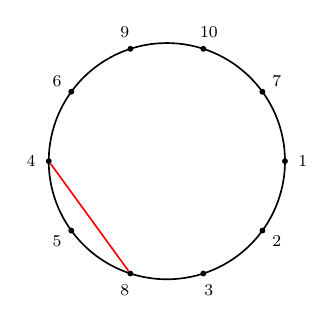}
  \end{subfigure}
  \caption{Example of sign patterns after each iteration in Algorithm~\ref{alg: main}}
  \label{fig:ExAlg}
\end{figure}
\end{example}

The proof that Algorithm~\ref{alg: main} terminates when $s$ is consistent goes roughly as follows. Let $N(s)$ denote the number of negative chords in $s$, and $\ell(s)$ the length of the shortest negative chord. 
The pair $(N, \ell)$ is an invariant assigned to each sign pattern $s$. Our strategy is to show that this invariant decreases.
First, we show that if $\ell(s)=2$ or $3$, then one iteration of the while loop produces a sign pattern $t$ with $N(t)<N(s)$. When $\ell(s)>3$ we show that after a finite number of iterations the algorithm produces a sign pattern $t$ with $N(t)\leq N(s)$ and $\ell(t)<\ell(s)$.
\smallskip

We split the remainder of this section into four parts. In Section~\ref{subsec: 3.1} we consider the case $\ell(s)=2$. In Section~\ref{subsec: map} we compute the monomial map induced by the transposition $(1 \ell)$. In Section~\ref{subsec: 3.2} we use this to show that, in most cases where $\ell \geq 3$, the invariant $(N, \ell)$ is decreased after one iteration of the while loop.
Section~\ref{subsec: 3.3} deals with the remaining cases.

\subsection{Negative chords of length two}\label{subsec: 3.1}

From now on, $s$ will be a sign pattern in the standard dihedral coordinates on $\M_{0,n}$. Without loss of generality we may assume that the shortest negative chord in $s$ is $n\ell$. In this section $\ell=2$.

\begin{proposition}\label{lem: reducing negatives when j = 2}
    Let $s$ be a consistent sign pattern such that $s(u_{n2}) = -$. Let $t$ be the sign pattern given by $\varphi_{(12)}(\M_{0,n}(\RR_s)) = \M_{0,n}^{(12)}(\RR_t)$. Then $N(t) < N(s)$.
\end{proposition}

We need the following technical result which guarantees that ``problematic'' sign patterns, i.e. the sign patterns $s$ for which $N(s) \leq  N(t)$, cannot be consistent.

\begin{lemma}\label{lem: technical lemma}
    Let $s$ be a sign pattern with $s(u_{n2})=-$. If any of the following holds
    \begin{enumerate}[label=(\roman*)]
        \item $s(u_{in}, u_{i1}, u_{i2}) = (+, -, +)$ for some $i \in [n]\setminus\{1,2,3,n-1,n\}$,
        \item $s(u_{1\,n-1}, u_{2\,n-1}) = (-,+)$,
        \item $s(u_{3n}, u_{13}) = (+,-)$,
    \end{enumerate}
    then $s$ is not consistent. 
\end{lemma}

\begin{figure}[h!]
  \centering
  \begin{subfigure}[c]{0.33\textwidth}
  	\includegraphics[width=\textwidth]{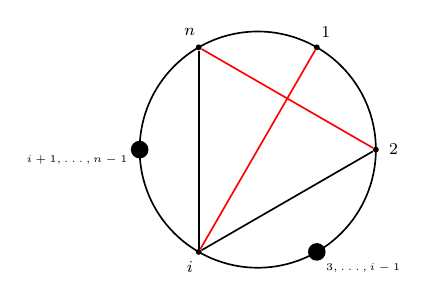}
    \caption{$(i)$}
    \label{fig:i}
  \end{subfigure}%
  \hfill %
  \begin{subfigure}[c]{0.26\textwidth}
  \centering
	\includegraphics[width=\textwidth]{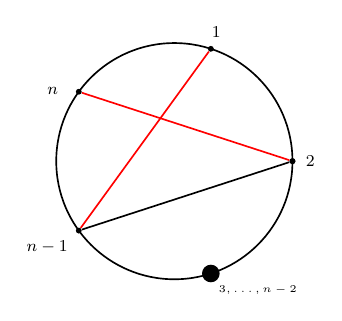}
    \caption{$(ii)$}
    \label{fig:ii}
  \end{subfigure}%
  \hfill%
  \begin{subfigure}[c]{0.28\textwidth}
  \centering
  	\includegraphics[width=\textwidth]{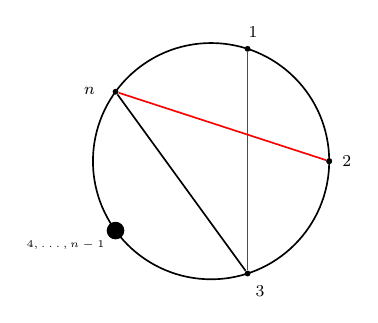}
    \caption{$(iii)$}
    \label{fig:iii}
  \end{subfigure}
  \caption{The cases of Lemma~\ref{lem: technical lemma}, where negative chords are shown in {\color{red}red} and positive chords are shown in black.}
  \label{fig: pentagonhexagon}
\end{figure}

\begin{proof}
We look for a collection of extended $u$-relations which contradict the assumptions in the lemma. 
Let us consider case $(ii)$ first.
It is sufficient to look for relations $R_{A,B,C,D}$ such that the set of adjacent indices $\{3, \ldots, n-2\}$ is contained in one of $A,B, C$ or $D$.
We regard these extended relations as primitive $u$-relations for $\M_{0,5}$, where the points are labeled $n, 1, 2, \{3, \ldots, n-2\}, n-1$, as shown in Figure~\ref{fig:ii}.
Here the chord $u_{I J}$ is the product $\prod_{i \in I, \, j\in J} u_{ij}$ of $u$-variables coming from the $n$-gon.
We check that $s$ cannot restrict to any of the 12 consistent sign patterns on $\M_{0,5}$.
The relation $R_{n2}$ implies that $s(u_{1 \, \{3 \ldots n-2\}}) = -$, because $s(u_{n2}) = - = s(u_{n-1 \, 1})$ by assumption. 
But this contradicts the relation $R_{1\, \{3\ldots n-2\}}$, since we also have $s(u_{n-1 \, \{1,2\}}) = -$ by our assumptions.
The argument that $(iii)$ leads to an inconsistency is symmetric to the argument for $(ii)$: this time, the vertices of the pentagon are labeled $n, 1, 2, 3, \{4, \ldots, n-1\}$ -- see Figure~\ref{fig:iii}.

For case $(i)$, we consider a hexagon rather than a pentagon: we label its vertices by $n, 1, 2, S := \{3, \ldots, i-1\}, i, T := \{i+1, \ldots, n-1\}$, as shown in Figure~\ref{fig:i}. 
One can deduce the signs of the chords on this hexagon by analyzing the sequence of relations $R_{S T}, R_{n2}, R_{1S}, R_{1T}, R_{1i}, R_{n | 12 | Si|T}$. The resulting sign pattern is inconsistent.
\end{proof}

\smallskip

\begin{proof}[Proof of Proposition~\ref{lem: reducing negatives when j = 2}]
    We compare the sign patterns $s$ and $t$ via our computation of the monomial map $\phi_{(12)} : \CC[\M_{0,n}^{(12)}] \longrightarrow \CC[\M_{0,n}]$ from Example~\ref{ex: monomial map for (12)}. 
    Explicitly, we use the fact that $t(v_{km}) = s(\phi_{(12)}(v_{km}))$, where we are continuing with the notation $v_{(12)(i)\,(12)(j)} = u_{ij}^{(12)}$ for the dihedral coordinates on $\M_{0,n}^{(12)}$ set in the example.
    
    For $i, j \in [n] \setminus \{n, 1 , 2\}$, we have $t(v_{ij}) = s(u_{ij})$.
    After applying the transposition $(12)$, the length two chord $u_{n2}$ becomes the length two chord $v_{n1}$ and we have
    $$
    t(v_{n1}) = s\left( - u_{n2} \prod_{j \notin \{1, 2, n\}} u_{2j}^{-1}\right) = - s(u_{n2}) s \left(\prod_{j \notin \{1, 2, n\}} u_{2j}\right) 
    $$
    By assumption $s(u_{n2}) = -$, and so the $u$-relation $R_{2n}$ guarantees that $s(\prod_{j \notin \{1, 2, n\}} u_{2j}) =+$. Thus $t(v_{n1}) = +$. 
    The number of negative chords in $t$ restricted to the set
    $\{v_{ij}: i,j \neq 1,2,n\} \cup \{v_{n1}\}$ is thus exactly one less than in $s$ restricted to $\{u_{ij}: i,j \neq 1,2,n\} \cup \{u_{n2}\}$. The remaining coordinates can be partitioned into pairs of bijective subsets
    \begin{align*}
        A_i := \{u_{in}, u_{1i}, u_{2i}\} &\longleftrightarrow \{v_{in}, v_{2i}, v_{1i}\} =: B_i \quad \text{for all } i \notin \{1,2,3,n-1,n\}\\
        A_{n-1} := \{u_{n-1\,1}, u_{n-1\,2}\} & \longleftrightarrow \{v_{n-1\,2}, v_{n-1\,1}\} := B_{n-1} \\
        A_{3} := \{u_{n3}, u_{13}\} & \longleftrightarrow \{v_{n3}, v_{23}\} := B_{3}.
    \end{align*}
    We finish the proof by showing that for each pair $(A_i, B_i)$, the number of negatives in $B_i$ is less than or equal to the number of negatives in $A_i$. For each $i \notin \{1,2,n\}$ we have 
    $$
    t(v_{in}) = s(u_{1i})s(u_{in}), \quad t(v_{1i}) = s(u_{1i}) s(u_{2i}) \quad \text{and} \quad t(v_{2i}) = s(u_{1i}).
    $$
    It follows that $B_i$ contains more negatives than $A_i$ if and only if 
    \begin{equation}\label{eq: problem case}
      s(u_{in}, u_{1i}, u_{2i}) = (+, -, +)         
    \end{equation}
    where these chords exist. By Lemma~\ref{lem: technical lemma}, \eqref{eq: problem case} contradicts the consistency of $s$.
\end{proof}

\subsection{The monomial map $\mathbf{\phi_{(1\ell)}}$}\label{subsec: map}
Now we let $s$ be a consistent sign pattern with shortest negative chord $u_{n\ell}$, where $\ell\geq 3$. In order to understand the sign pattern obtained after the first step of Algorithm~\ref{alg: main}, we compute the isomorphism $\phi_{(1\ell)}: \CC[\M_{0,n}^{(1\ell)}] \rightarrow \CC[\M_{0,n}]$ (see the proof of Lemma~\ref{lem:monomial_map} for the definition). 
For ease of notation, we write the coordinates on $\M_{0,n}^{(1\ell)}$ as $v_{ij}$ -- that is, $v_{ij}=v_{(1\ell)(k)\, (1\ell)(p)} := u_{kp}^{(1\ell)}$ as defined in~\eqref{eq: notation uij}.

\begin{lemma}\label{lem: maps-formulas}
    Table~\ref{tab:technical_comp} gives the image of each variable $v_{ij}$ under the monomial map $\phi_{(1\ell)}$.
    \afterpage{
    \vspace*{\fill}
    \begin{center}
    \begin{table*}[h]
        \centering
        \renewcommand{\arraystretch}{2.5}
        \begin{tabular}{ccc|c}
             \multirow{2}{*}{$v_{ij}$} & \multirow{2}{*}{$\longmapsto$}  & \multirow{2}{*}{$u_{ij}$} &$i,j \in \{2, \ldots, \ell-2\}$  or\\
            &&& $\begin{array}{l}
                i \in \{\ell+1, \ldots, n-1\}, \vspace{-1.7em}\\
                j \notin \{n,1,\ell-1,\ell\}
            \end{array}$
             \\
            \hline
             $v_{jn}$ & $\longmapsto$ & $\prod\limits_{k = \ell}^{n-1} u_{kj}^{-1}$ & \multirow{2}{*}{$j \in \{2, \ldots, \ell-2\}$ }\\
            $v_{1j}$ & $\longmapsto$ & $\prod\limits_{k = \ell+1}^n u_{kj}^{-1}$ &\\
            \hline
             $v_{j\ell}$ & $\longmapsto$ & $u_{1j} \prod\limits_{k = \ell}^{n} u_{kj}$ & $j \in \{3, \ldots, \ell-2\}$ \\
            \hline
             $v_{j \, \ell-1}$ & $\longmapsto$ & $\prod\limits_{k=\ell-1}^n u_{kj}$ & $j \in \{2, \ldots, \ell-3\}$\\
            \hline
            $v_{1i}$ & $\longmapsto$ & $\prod\limits_{j=1}^\ell u_{ij}$ & \multirow{4}{*}{$i \in \{\ell+1, \ldots, n-1\}$}\\     $v_{\ell-1 \, i}$ & $\longmapsto$ & $\prod\limits_{j=1}^{\ell-2} u_{ij}$ &\\  
            $v_{\ell i}$ & $\longmapsto$ & $\prod\limits_{j=2}^{\ell-1} u_{ij}^{-1}$ &\\  
            $v_{in}$ & $\longmapsto$ & $ \prod\limits_{j = n}^{\ell-1} u_{ij}$ &\\  
            \hline
            $v_{\ell1}$& $\longmapsto$ &$\prod\limits_{kh \not \sim 1\ell} u_{kj}$ &\\
            $v_{n1}$& $\longmapsto$ &$- u_{n\ell} \prod\limits_{kj \not \sim n\ell} u_{kj}^{-1}$ &\\
            $v_{\ell-1 \, \ell}$& $\longmapsto$ & $-  u_{1 \, \ell-1} \prod\limits_{kj \not \sim 1 \, \ell-1} u_{kj}^{-1}$ &\\
            $v_{\ell-1 \, n}$& $\longmapsto$ &$\prod\limits_{kj \not \sim n \, \ell-1} u_{kj}$ &
        \end{tabular}
        \renewcommand{\arraystretch}{1}
        \caption{The images of the monomial map $\phi_{(1\ell)}$}\label{tab:technical_comp}
    \end{table*}
    \end{center}
    \vspace*{\fill}
    \clearpage
}
\end{lemma}

\begin{proof}
    The transposition $(1\ell)$ applied to the standard dihedral ordering can be achieved by the following sequence of $L:=2\ell-3$ adjacent transpositions:
    \begin{align}\label{eq: 1j as adj}
        (1\ell) & = (2\ell) (3\ell) \cdots (\ell-1 \, \ell) (1\ell) (1 \, \ell-1) \cdots (13) (12)\\
        & = \tau_1^{(2\ell)(1\ell)} \circ \tau_2^{(3\ell)(2\ell)(1\ell)} \circ \cdots \circ \tau_{\ell-2}^{(1\ell) \cdots (12)} \circ \tau_{\ell-1}^{(1 \, \ell-1) \cdots (12)} \circ \tau_{\ell-2}^{(1 \, \ell-2) \cdots (12)}\circ \cdots \circ \tau_2^{(12)} \circ \tau_1. \nonumber
    \end{align}
    Recall that to construct the map $\phi_{(1\ell)}$ we ``sort" the dihedral ordering represented by $(1\ell) \in S_n$. This means applying the transpositions in~\eqref{eq: 1j as adj} from left to right to $(1\ell)$ to recover $\text{id} \in S_n$.
    
    For each $m \in [L]$, let $\sigma_m$ denote the composition of the first $m$ adjacent transpositions in this sequence from right to left -- i.e. $\sigma_1 = (12), \sigma_2 = (13)(12), \ldots$ and so on until $\sigma_{L} = (1\ell)$. The monomial transformation $\CC[\M_{0,n}^{(1\ell)}] \rightarrow \CC[\M_{0,n}]$ is given by the composition
    $$
    \phi_{(1\ell)} = \phi_1^{\text{id}} \circ \phi_2^{\sigma_1} \circ \cdots \phi_{\ell-2}^{\sigma_{\ell-3}} \circ \phi_{\ell-1}^{\sigma_{\ell-2}} \circ \phi_{\ell-2}^{\sigma_{\ell-1}} \circ \cdots \circ \phi_3^{\sigma_{L-2}} \circ \phi_2^{\sigma_{L-1}}.
    $$
    One computes $\phi_{(1\ell)}(v_{ij})$ by successively applying the formula for $\phi_k^{\sigma}$ (see Definition~\ref{def: monomial map for general adjacent transposition}).
    As an example, we give the step-by-step computation of $\phi_{(1\ell)}(v_{\ell-1 \, \ell})= -  u_{1 \, \ell-1} \prod_{kj \not \sim 1 \, \ell-1} u_{kj}^{-1}$ below. All the other entries in Table~\ref{tab:technical_comp} are computed similarly. 
    \begin{align*}
         &v_{\ell -1\, \ell} \; \xmapsto{\phi_{2}^{\sigma_{L-1}}} \; (u^{\sigma_{L-1}}_{1\,\ell-1})^{-1} \; \xmapsto{\phi_{3}^{\sigma_{L-2}}} \; (u_{1\,\ell-1}^{\sigma_{L-2}}u_{2\,\ell-1}^{\sigma_{L-2}})^{-1}\; \xmapsto{\phi_{4}^{\sigma_{L-3}}} \; \left(u_{1\,\ell-1}^{\sigma_{L-3}} u_{2\,\ell-1}^{\sigma_{L-3}} u_{3\,\ell-1}^{\sigma_{L-3}}\right)^{-1} \; \ldots \;\\ 
         &\xmapsto{\phi_{\ell-3}^{\sigma_{\ell}}} \; \prod_{k = 1}^{\ell-3} (u_{k\,\ell-1}^{\sigma_{\ell}})^{-1} = \left(u_{\ell-3\,\ell-1}^{\sigma_{\ell}} \prod_{k = 1}^{\ell-4} u_{k\, \ell-1}^{\sigma_{\ell}}\right)^{-1}\\
        & \xmapsto{\phi_{\ell-2}^{\sigma_{\ell-1}}} \; - \; (u_{\ell-3\, \ell-1}^{\sigma_{\ell-1}})^{-1} \prod_{k \notin \{\ell-3, \ell-2, \ell-1\}}u_{k\,
        \ell-2}^{\sigma_{\ell-1}} \prod_{k = 1}^{\ell-4} (u_{k\,\ell-1}^{\sigma_{\ell-1}} u_{k\,\ell-2}^{\sigma_{\ell-1}})^{-1} = \; -\prod_{k=1}^{\ell-3} (u_{k\,\ell-1}^{\sigma_{\ell-1}})^{-1}\prod_{k=\ell+1}^{n} u_{k\,\ell-2}^{\sigma_{\ell-1}} u_{\ell\,\ell-2}\\
        &\xmapsto{\phi_{\ell-1}^{\sigma_{\ell-2}}} \; -  \prod_{k=1}^{\ell-3} u_{k\,\ell-1}^{\sigma_{\ell-2}} \prod_{k=\ell+1}^{n} u_{k\,\ell-2} u_{k\,\ell-1}^{\sigma_{\ell-2}} (-u_{\ell\,\ell-2}^{\sigma_{\ell-2}}) \prod_{k \notin \{\ell-2, \ell-1, \ell\}} (u_{k\,\ell-1}^{\sigma_{\ell-2}})^{-1} = \; \prod_{k = \ell}^n u_{k\,\ell-2}^{\sigma_{\ell-2}} \\ 
        &\xmapsto{\phi_{\ell-2}^{\sigma_{\ell-3}}} \; \prod_{k = \ell}^n (u_{k\,\ell-2}^{\sigma_{\ell-3}})^{-1} \; \xmapsto{\phi_{\ell-4}^{\sigma_{\ell-4}}} \; \prod_{k = \ell}^n (u_{k\,\ell-3}^{\sigma_{\ell-4}} u_{k\,\ell-2}^{\sigma_{\ell-4}})^{-1}\\  
        &\vdots \\
        &\xmapsto{\; \phi_{2}^{\sigma_{1}}\; } \; \prod_{k = \ell}^n (u_{k2}^{\sigma_{1}} \cdots u_{k\,\ell-2}^{\sigma_{1}})^{-1} = \left(u_{n2}^{\sigma_{1}}u_{n3}^{\sigma_{1}} \cdots u_{n\,\ell-2}^{\sigma_{1}}\prod_{k = \ell}^{n-1} u_{k2}^{\sigma_{1}} \prod_{k = \ell}^{n-1} u_{k3}^{\sigma_{1}} \ldots u_{k\,\ell-2}^{\sigma_{1}}\right)^{-1} \\
        &\xmapsto{\; \phi_{1}^{\text{id}}\;} \;  - \, u_{n2}^{-1} \prod_{k \notin \{n,1,2\}} u_{k1} \prod_{k=3}^{\ell-2} u_{kn}^{-1}u_{k1}^{-1}
        \prod_{k=\ell}^{n-1} u_{k1} \prod_{k=\ell}^{n-1} \prod_{j = 3}^{\ell-2} u_{kj}^{-1} = - u_{1 \, \ell-1} \prod_{k=\ell}^{n} \prod_{j=2}^{\ell-2} u_{kj}^{-1}. \qedhere
    \end{align*}
\end{proof}

\subsection{Decreasing the invariant}\label{subsec: 3.2}

We now use the monomial map given in Table~\ref{tab:technical_comp} to sign chase from $\M_{0,n}$ to $\M_{0,n}^{(1\ell)}$. 
To show that the invariant decreases, we deploy the same strategy as in Section~\ref{subsec: 3.1} and partition the two sets of variables $\{v_{ij}\}$ and $\{u_{ij}\}$ into pairs of bijective subsets.
We compare the signs of the $v$'s and $u$'s in each pair separately. 

\medskip

Recall that $s$ is consistent sign pattern in the $u$'s such that the shortest negative chord is $u_{n\ell}$. Let $t$ be the preimage of $s$ under the monomial map $\phi_{(1\ell)}$, or equivalently, the sign pattern of $\varphi_{(1\ell)}(\M_{0,n}(\RR_s))$.
The first pair of subsets we consider comes from the first row of Table~\ref{tab:technical_comp}. It is clear that $s$ and $t$ have the same number of negatives among the two sets of variables below.
\begin{equation*}
    \left\{u_{ij} \,\middle| \begin{array}{l}
        i,j \in \{2, \ldots, \ell-2\}\}, \text{or}\\
         i \in \{\ell+1, \ldots, n-1\}\\
         \text{and } j \notin \{n,1,\ell-1,\ell\}
    \end{array}  \right\}   \longleftrightarrow   \left\{v_{ij} \, \middle| \begin{array}{l}
        i,j \in \{2, \ldots, \ell-2\}\}, \text{or}\\
         i \in \{\ell+1, \ldots, n-1\}\\
         \text{and } j \notin \{n,1,\ell-1,\ell\}
    \end{array}\right\}
\end{equation*}
The next collection of pairs we consider are
\begin{align}
    \{u_{n\ell}, u_{n \, \ell-1}, u_{1\ell}, u_{1\, \ell-1}\}  & \longleftrightarrow  \{v_{n1}, v_{n \, \ell-1}, v_{\ell1}, v_{\ell \, \ell-1}\} \label{eq: partitionpair1}\\
    \{u_{nj}, u_{1j}, u_{j \, \ell-1 }, u_{j \ell}\}  & \longleftrightarrow   \{v_{nj}, v_{\ell j}, v_{j \, \ell-1}, v_{j1}\} \quad \text{for } j \in \{2, \ldots, \ell-2\}. \nonumber
\end{align}
The following lemma shows that $s$ and $t$ also have the same number of negatives among these pairs of subsets.

\begin{lemma}\label{lem: pos_in_1j}
    The sign pattern $t$ takes the following values
$$
t(v_{n1}, v_{n \, \ell-1}, v_{\ell1}, v_{\ell \, \ell-1}) = (+, +, +, -) \qquad \text{and} \qquad  t(v_{nj}, v_{\ell j}, v_{j \, \ell-1}, v_{j1}) = (+, +, +, +)
$$
for all $j=2,\ldots,\ell-2$. Here we consider $v_{\ell2}$ and $v_{\ell-2 \, \ell-1}$ to be positive.
\end{lemma}

\begin{proof}
The signs of $t$ on the $v$-coordinates can be read off from Table~\ref{tab:technical_comp}, since $t(v_{ij}) = s(\phi_{(1\ell)}(v_{ij}))$.
We will use the following shorthand notations 
\begin{equation}\label{eq: monomial Mi}
    M_i := \prod_{k = 2}^{\ell-2} u_{ik}, \quad N_j := \prod_{k = \ell+1}^{n-1} u_{kj} \quad \text{and} \quad L := \prod_{i = \ell+1}^{n-1} M_i = \prod_{j = 2}^{\ell-2} N_j.
\end{equation}
Recall that $s$ is positive on all chords of length strictly less than $\ell$, so
$$
s(u_{n \, \ell-1}) = s(u_{1\, \ell-1}) = s(u_{1\ell}) = s(M_n) = s(M_\ell) = +.
$$
Since $s(u_{n\ell}) = -$, it follows from the $u$-relation $R_{n\ell}$ that $s(\prod_{kj \not \sim n\ell} u_{kj}) = +$. Thus,
$$
t(v_{n1}) = s\left(- u_{n\ell} \prod\limits_{kj \not \sim n\ell} u_{kj}^{-1}\right) = - s(u_{n\ell})  s\left(\prod_{kj \not \sim n\ell} u_{kj}\right) = +.
$$ 
The extended $u$-relation
$
(u_{n\, \ell-1}u_{n\ell} + N_1 L - 1)
$
implies that $s(N_1 L) = +$, and therefore
$$
t(v_{n \, \ell-1}) = s(u_{1\ell}M_\ell N_1L) = +.
$$
Similarly, the relation $(u_{n\ell}u_{1\ell} + N_{\ell-1}L - 1)$ implies that 
$
t(v_{\ell1}) = s(u_{n \, \ell-1} M_n N_{\ell-1} L) = +,
$
and $(u_{n \, \ell-1} u_{n\ell} u_{1 \, \ell-1} u_{1\ell} + L - 1)$ implies that 
$
t(v_{\ell \, \ell-1}) = s(-u_{1\, \ell-1} M_n L M_\ell) = -.
$

\medskip
Now we determine the signs $t(v_{nj}, v_{\ell j}, v_{j \, \ell-1}, v_{j1})$.  
Note that the coordinates $u_{nj}, u_{1j}$, $u_{\ell-1 \, j}$ and $u_{\ell j}$ are associated to chords of length less than $\ell$ when $j \in \{2, \ldots, \ell-2\}$. 
Furthermore, $s(N_j) = +$ for all $j \in \{2, \ldots, \ell-2\}$. 
Indeed, $s$ is negative on the monomial
$$
\prod_{k = n}^{j-1} \prod_{i = j+1}^\ell u_{ik},
$$
since $u_{n\ell}$ is the only factor associated to a chord of length $\ell$ or greater. 
The extended $u$-relation $R_{n, \ldots, j-1 | l | j+1, \ldots, \ell | \ell+1, \ldots, n-1}$ therefore implies that $s(N_j) = +$.
Thus we have
$$
t(v_{nj}, v_{\ell j}, v_{j \, \ell-1}, v_{j1}) = s\left(u_{\ell j}N_j, u_{nj}u_{1j}u_{j \ell} N_j, u_{nj}u_{j\,\ell-1}u_{l\ell}N_j, u_{nj}N_j\right) = (+,+,+,+). \qedhere
$$
\end{proof}

\smallskip

The remaining variables can be partitioned into these pairs of subsets
$$
\{u_{ni}, u_{1i}, u_{\ell-1 \, i}, u_{\ell i}\} \longleftrightarrow \{v_{ni}, v_{\ell i}, v_{\ell-1 \, i}, v_{1i}\} \quad \text{for $i \in \{\ell+1, \ldots, n-1\}$.}
$$
We see from Table~\ref{tab:technical_comp} that the signs of the four coordinates $v_{in}, v_{1i}, v_{\ell-1 \, i}$ and $v_{\ell i}$ depend not only on the signs of $u_{in}, u_{1i}, u_{\ell-1 \, i}$ and $u_{\ell i}$, but also on the sign of the monomial $M_i$~\eqref{eq: monomial Mi}.
By considering all possible sign combinations of $M_i, u_{in}, u_{1i}, u_{\ell-1 \, i}$ and $u_{\ell i}$ we see that the only ones which increase the number of negative chords, i.e. $N(t)>N(s)$, are shown in Table~\ref{tab:mult table Mi}.
\begin{table}[h!]
    \centering
    \begin{tabular}{|ccccc|cccc|c|}
        \hline
         $M_i$ & $u_{in}$ & $u_{1i}$ & $u_{\ell-1\,i}$ & $u_{\ell i}$ & $v_{in}$ & $v_{1i}$ & $v_{\ell-1\,i}$ & $v_{\ell i}$ & \\
         \hline
        $-$ & + & + & + &+ & $-$ &$-$ & $-$ &$-$ & \text{not consistent}\\
        + & + & $-$ & + &+ & $-$ &$-$ & $-$ &+ & \text{not consistent}\\
        + & + & + & $-$ &+ & $-$ &$-$ & + &$-$ & \text{not consistent}\\
        $-$ & $-$ & + & + &+ & + &$-$ & $-$ &$-$ & \text{consistent}\\
        $-$ & + & + & + &$-$ & $-$ &+ & $-$ &$-$ & \text{consistent}\\
        \hline
    \end{tabular}
    \caption{Problematic sign patterns for $(1\ell)$.}
    \label{tab:mult table Mi}
\end{table}

We have proved the following generalization of Proposition~\ref{lem: reducing negatives when j = 2} to the case $\ell>2$.
\begin{proposition}\label{lem: reducing negatives for general j}
If $s$ avoids the sign patterns listed in Table~\ref{tab:mult table Mi} for all $i \in \{\ell+1, \ldots, n-1\}$, then $N(t) \leq N(s)$ and $t$ has a shorter negative chord $v_{\ell \, \ell-1}$.
Thus the invariant has decreased.
\end{proposition}

\subsection{Problematic sign patterns}\label{subsec: 3.3} 
To complete the proof of Theorem~\ref{thm:mainResult}, we deal with the remaining cases listed in Table~\ref{tab:mult table Mi}. First we show that, under our assumptions, we may discount the majority of these problematic sign patterns.

\begin{lemma}\label{lem:incos}
    Let $s$ be a consistent sign pattern with shortest negative chord $u_{n\ell}$. Then $s$ cannot restrict to any of the first three rows in Table~\ref{tab:mult table Mi}.
\end{lemma}
\begin{proof}
    For $i \in \{\ell+2, \ldots, n-2\}$ we consider an octagon labeled by a partition of $[n]$ as in Figure~\ref{fig:8partition}, grouping together some of the vertices.
\begin{figure}[h!]
  \centering
  	\includegraphics[width=0.5\textwidth]{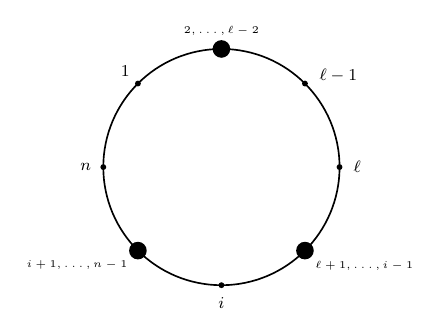}
    \caption{Partition of $n$ vertices into eight subsets}
    \label{fig:8partition}
\end{figure}
    It is enough to show that, if $s$ restricts to one of the first three rows of Table~\ref{tab:mult table Mi}, then it is not consistent on the octagon. Indeed, any extended $u$-relation for this octagon can be translated to an extended $u$-relation in the bigger $n$-gon. We use Macaulay2 - the code is accessible in \cite{M2Code}. Concretely, we compute all consistent sign patterns in the octagon and verify that none of them restricts to any of the first three rows of Table~\ref{tab:mult table Mi}.

    In the cases $i=\ell+1$ or $i=n-1$, we can rule out the first row: under the assumption that no chord of length less than $\ell$ is negative, $M_{\ell +1}$ and $M_{n-1}$ must be positive.
    The other two rows become
    \[
    s(M_{\ell+1}, u_{n \, \ell+1}, u_{1\,\ell+1}, u_{\ell-1 \, \ell+1}) = (+,+,-,+) \quad s(M_{n-1}, u_{1 \, n-1}, u_{\ell-1 \, n-1}, u_{\ell \, n-1}) = (+,+,-,+)
    \]
    because $u_{\ell \, \ell+1}$ and $u_{n-1 \, n}$ are not chords.
     These also contradict consistency. We verify this with Macaulay2 as above, using a heptagon instead of an octagon.
\end{proof}

\begin{corollary}[$\ell=3$]\label{rem:j=3}
Let $s$ be a consistent sign pattern with shortest negative chord $u_{n3}$. Let $t$ be the sign pattern given by $\varphi_{(13)}(\M_{0,n}(\RR_s)) = \M_{0,n}^{(13)}(\RR_t)$. Then $N(t) < N(s)$.
\end{corollary}

\begin{proof}
    In this case, the number of negatives in the sets of equation~\eqref{eq: partitionpair1} decreases by one as $v_{\ell \, \ell-1}$ is not a coordinate. Also, $M_i=1$ as there are not enough vertices between $1$ and $3$. Hence, all the problematic sign patterns which could increase the number of negatives are inconsistent by Lemma~\ref{lem:incos}.
\end{proof}

For the two last rows of Table~\ref{tab:mult table Mi}, where $M_i$ and another chord are negative, applying $(1\ell)$ actually can increase the number of negatives. However, this can be solved by applying more permutations. 
We will show that, after a finite number of iterations of the while loop in Algorithm~\ref{alg: main}, we end up with a sign pattern with at most as many negatives as $s$, and with a negative chord of length less than $\ell$.
In other words: we do not decrease the invariant in one iteration of the while loop, but after multiple iterations we do. 

\begin{lemma}\label{lem: final step}
    Let $s$ be a consistent sign pattern with shortest negative chord $u_{n\ell}$.
    Denote by $I \subset \{\ell+2, \ldots, n-2\}$ the set of indices $i$ for which $s$ restricts to one of the two last rows of Table~\ref{tab:mult table Mi}. We define 
    \[ k_i := \min\{1\leq p\leq \lfloor (\ell-3)/2 \rfloor: s(u_{p+1 \, i}) = - \text{ or }  s(u_{\ell-1-p \, i}) = -\} \quad \text{and} \quad k := \max_{i \in I}k_i. \]
    Let $\alpha := (k+1\,\ell-k)\cdots (3\,\ell-2)(2\,\ell-1) (1\,\ell)$ and $t$ be the preimage of $s$ under $\phi_\alpha$. 
    Then either $N(t) < N(s)$, or $N(t) = N(s)$ and $\ell(t) < \ell(s)$.
\end{lemma}

Before giving the proof, we illustrate the argument in the following example.

\begin{example}[$n=10$]\label{ex: n10}
    We consider the following consistent sign pattern in a decagon, shown in Figure~\ref{fig: decagon sign pattern}: the five colored (black, \textcolor{dgreen}{green}, \textcolor{dviolet}{violet}, \textcolor{dorange}{orange} and \textcolor{blue}{blue}) chords are negative and all other chords are positive. 
    \begin{figure}[h!]
        \centering
        \includegraphics[width=0.33\linewidth]{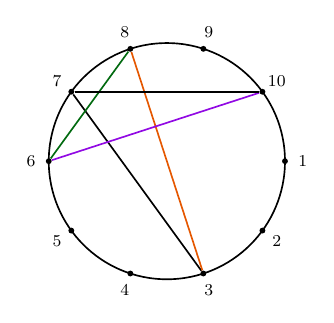}
        \caption{A consistent sign pattern in $n=10$}
        \label{fig: decagon sign pattern}
    \end{figure}
    We consider the negative chord $6\, 10$ - even though it is not the shortest negative chord, we can still apply the strategy of Lemma~\ref{lem: final step} to see what happens in the problematic cases. 
    Explicitly, we have $l=6$, $i = 8$, and
    \[ M_i = u_{28}u_{38}u_{48} < 0, \quad u_{68} <0 \quad \text{and} \quad u_{18},u_{58},u_{8\,10} >0, \]
    placing us in the last row of Table~\ref{tab:mult table Mi}.
    In the language of Lemma~\ref{lem: final step},  $k=2$ and $\alpha= (34)(25)(16)$.
    After applying each transposition sequentially we obtain the sign patterns in Figure~\ref{fig:ExMi}. The color coding helps to keep track of how much we have increased, and then decreased, the number of negative chords.
    \begin{figure}[h!]
  \centering
  \begin{subfigure}[c]{0.33\textwidth}
  	\includegraphics[width=\textwidth]{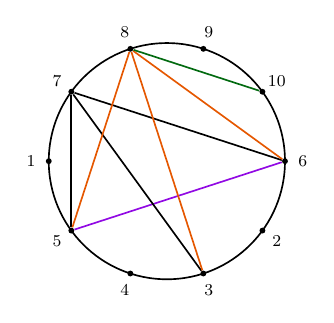}
    \caption{$(16)$}
    \label{fig:ExMi1}
  \end{subfigure}%
  \hfill %
  \begin{subfigure}[c]{0.33\textwidth}
  \centering
	\includegraphics[width=\textwidth]{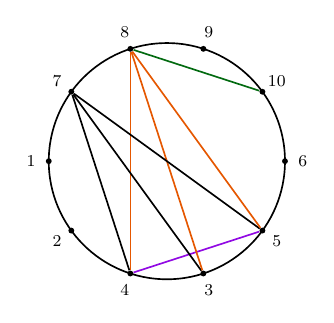}
    \caption{$(25)(16)$}
    \label{fig:ExMi2}
  \end{subfigure}%
  \hfill%
  \begin{subfigure}[c]{0.33\textwidth}
  \centering
  	\includegraphics[width=\textwidth]{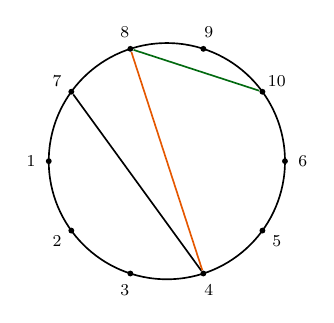}
    \caption{$(34)(25)(16)$}
    \label{fig:ExMi3}
  \end{subfigure}
  \caption{Sign patterns after iteratively applying transpositions.}
  \label{fig:ExMi}
\end{figure}
As expected from Table~\ref{tab:mult table Mi}, the number of negatives in the sign pattern~\ref{fig:ExMi1} has increased from $5$ to $8$. Applying the transposition $(2 5)$ does not change the number of negatives, but does shorten the negative chord $65$ to $54$, as shown in picture~\ref{fig:ExMi2}. Finally, applying the transposition $(3 4)$ does reduce the number of negative chords - see picture~\ref{fig:ExMi3}.
\end{example}

\medskip

\begin{proof}[Proof of Lemma~\ref{lem: final step}]
    In the two last rows of Table~\ref{tab:mult table Mi}, the number of negative chords increases in the following way: one negative chord in the original coordinates, $u_{in}$ or $u_{\ell i}$, corresponds to three negative chords in the new coordinates, $v_{1i},v_{\ell-1 \, i}, v_{\ell i}$ or $v_{in},v_{\ell-1 \, i}, v_{\ell i}$. 
    The two cases are completely symmetrical, so we only consider the first, that is $s(u_{in})=-$. 
    
    We keep track of the number of negative chords as we apply each transposition in $\alpha = (k+1\,\ell-k)\cdots (3\,\ell-2)(2\,\ell-1) (1 \ell)$. Note that, up to relabeling, we can use the computations from Lemma~\ref{lem: maps-formulas}. Let $\sigma = (2 \; \ell-1)(1\ell)$, so that the coordinates on $\M_{0,n}^{(2\,\ell-1)(1\ell)}$ are denoted by $u_{ij}^\sigma$. Let $r$ be the preimage of $s$ under $\phi_{(1\ell)}$ and $t$ the preimage of $r$ under $\phi_{(2 \, \ell-1)}^{(1 \ell)}$. Denoting $M_i' := \prod_{j=3}^{\ell -3} v_{ij}$, we have $r(M_i' v_{2i}v_{\ell-2 \, i}) = s(\phi_{(1\ell)}(M_i' v_{2i}v_{\ell-2 \, i}) ) = s(M_i) = -$. Table~\ref{tab:mult table z} shows the possible restrictions of $r$ to $(M_i', v_{\ell i}, v_{2i}, v_{\ell-2 \, i}, v_{\ell-1 \, i})$, and the consequent restriction of $t$.
    \begin{table}[h!]
    \centering
    \begin{tabular}{|ccccc|cccc|}
        \hline
         $M_i'$ & $v_{\ell i}$ & $v_{2i}$ & $v_{\ell-2\,i}$ & $v_{\ell-1\,i}$ & $u_{1i}^\sigma$ & $u_{\ell-1\,i}^\sigma$ & $u_{\ell-2\,i}^\sigma$ & $u_{2i}^\sigma$\\
         \hline
         $-$ & $-$ & $-$ & $-$ &$-$ & + &+& + &+ \\
        + & $-$ & + & $-$ &$-$ & + &+ & +&$-$ \\
        + & $-$ & $-$ & + &$-$ & + &+ & $-$ &+ \\
        $-$ & $-$ & + & + &$-$ & + &+ & $-$ &$-$ \\
        \hline
    \end{tabular}
    \caption{Sign patterns for $(2\,\ell-1)(1\ell)$}
    \label{tab:mult table z}
\end{table}

In all rows of Table~\ref{tab:mult table z} except the last, we lose at least as many negatives as were added by $(1\ell)$.
In other words, in the first three rows of Table~\ref{tab:mult table z} we have $N(t)\leq N(s)$.
Applying Section~\ref{subsec: 3.2} to the transposition $(2 \, \ell-1)$ guarantees that, away from the variables listed in Table~\ref{tab:mult table z}, the number of negatives remains constant from $r$ to $t$.

\smallskip
We now consider the last row of Table~\ref{tab:mult table z}. In this case we have $N(t)=N(r)$ and $t(u_{\ell-2\,i}^\sigma)=t( u_{2i}^\sigma) = -$. 
This ensures that we can repeat a similar argument for the next transposition in $\alpha$, $(3\,\ell-2)$, and obtain the same sign table with columns labeled
$$
M_i''\; u^\sigma_{\ell-1 i}\; u^\sigma_{3i}\; u^\sigma_{\ell-3 \, i}\; u^\sigma_{\ell-2 \, i} \; | \; u_{2i}^{(3 \, \ell -2)\sigma}\; u_{\ell-2 \, i}^{(3 \, \ell -2)\sigma}\, u_{\ell-3 \, i}^{(3 \, \ell -2)\sigma}\; u_{3i}^{(3 \, \ell -2)\sigma}.
$$ 
We can reapply the results from Section~\ref{subsec: 3.2} to the sign pattern $t$ even though $u_{2 \, \ell-2}^\sigma$ may not be its shortest negative chord. In fact, if we consider the subdivision of the $n$-gon into two smaller polygons given by splitting along chord $u_{2 \, \ell-2}^\sigma$, we only require all chords in the smaller polygon to be positive.

For an index $i$ such that $k_i < k$, applying the transpositions $((k_i+1)+1 \; \ell -(k_i+1)),\ldots,$ $(k+1 \, \ell-k)$ does not increase the number of negative chords with vertex $i$. This is because the right column of Table~\ref{tab:mult table z} does not follow any of the problematic sign patterns from Table~\ref{tab:mult table Mi}. Thus, the definition of $k$ guarantees that after applying $(k+1\,\ell-k)$ the number of negative chords is at most $N(s)$.
Moreover, the shortest negative chord has length strictly smaller than $\ell(s)$, unless $u_{k \, \ell -k}^{\beta}$ has length two or three, where $\beta = (k\,\ell-k+1)\cdots (3\,\ell-2)(2\,\ell-1) (1\,\ell)$.
In this case, we end up with strictly fewer than $N(s)$ negative chords by Proposition~\ref{lem: reducing negatives when j = 2} and Corollary~\ref{rem:j=3}.
\end{proof}

We put all of the results together to complete the proof of our main theorem.

\begin{proof}[Proof of Theorem~\ref{thm:mainResult}]
    We prove that Algorithm~\ref{alg: main} terminates for a consistent sign pattern $s$.
    We proceed by double induction on $\ell = \ell(s)$ and $N = N(s)$ -- it suffices to show that finitely many iterations of the while loop will decrease one of them while not increasing the other.
    For $\ell= 2$ or $3$, we showed in Proposition~\ref{lem: reducing negatives when j = 2} and Corollary~\ref{rem:j=3} that $N$ always decreases after one iteration.
    For $\ell > 3$, we use Proposition~\ref{lem: reducing negatives for general j} and Lemma~\ref{lem: final step}. 

    We claim that the treatment of the problematic sign patterns, given in Lemma~\ref{lem: final step} follows the procedure prescribed by Algorithm~\ref{alg: main}. This is not immediately obvious -- after the transposition $(1\ell)$ has been applied, the negative chord $u^{(1 \ell)}_{1 \, \ell-1}$ might not be the shortest. 
    The only other new negative chords that could be created by the map $\phi_{(1\ell)}$ are of the form $u^{(1 \ell)}_{\ell \, i}$, $u^{(1 \ell)}_{in}$, $u^{(1 \ell)}_{\ell-1 \, i}$ and $u^{(1 \ell)}_{1 i}$ for $i \in \{\ell+2, \ldots, n-2\}$, the shorter of which are $u^{(1 \ell)}_{\ell \, i}$ and $u^{(1 \ell)}_{in}$.
    Then the next step in Algorithm~\ref{alg: main} is applying the one of the transpositions $(\ell+1 \, i)$ or $(i+1 \, n)$. 
    Both of these commute with the rest of the transpositions in $\alpha$ from Lemma~\ref{lem: final step}, i.e. they commute with the permutation $\sigma=(k+1\,\ell-k)\cdots (3\,\ell-2)(2\,\ell-1)$.
    Similarly, if $u_{2 \, \ell-2}^{(2\, \ell-1)(1 \ell)}$ is not a shortest chord in the sign pattern obtained after applying $(2 \, \ell-1)$, then the shortest negative must be of the form $u_{\ell-1 i}^{(2 \, \ell -1)(1 \ell)}$ or $u_{i1}^{(2 \, \ell -1)(1 \ell)}$ for some $i \in \{\ell +1, \ldots, n-1\}$.
    The next transposition applied in Algorithm~\ref{alg: main} is then $(1 i)$ or $(i+1 \, \ell)$, which both commute with $(k+1 \; \ell-k) \cdots (3 \, \ell-2)$.
    Arguing in the same manner for each transposition in the sequence $\alpha$ proves our claim.

    We have now proven that Algorithm~\ref{alg: main} terminates for a consistent sign pattern -- this proves that every consistent sign pattern is realizable.
    Proposition~\ref{prop: section 2} implies the rest of the statement of Theorem~\ref{thm:mainResult}.
\end{proof}

\section{Further directions}\label{sec:4}
In this section we discuss two further directions and related conjectures which appear naturally from our results. One is combinatorial as it extends our computations to a generalization of $\M_{0,n}$ using cluster algebras, see \cite{ClusterConfigSpaces}; while the other is algebro-geometric and asks for the prime ideal of the $u$-relations in the polynomial ring.

\subsection{Type C cluster configuration space} 
We start with a $2(n-1)$-gon whose vertices are labeled in pairs as $1,\ldots,n-1,\overline{1},\ldots,\overline{n-1}$. The type $C_n$ cluster configuration space $\M_{C_n}$ is given by folding from $\M_{0,2(n-1)-3}$. Concretely, its dihedral coordinates are indexed by pairs of chords $[i,j]=((i,j),(\overline{i},\overline{j}))$, with $i \in \{1,\ldots n-1\}$ and $j \in \{1,\ldots \overline{n-1}\}$, or centrally symmetric chords $[i,\overline{i}]$. The (extended) $u$-relations come from replacing the $u$-variables for $\M_{0,2(n-1)-3}$ by their corresponding $u$-variable in type $C$, i.e. labeled by the equivalence class of the chord. As this configuration space is closely related to the one studied in this article, also see \cite[Section 11.3]{ClusterConfigSpaces}, we believe our results should extend.

\begin{conjecture}
    The number of consistent sign patterns for the type $C$ configuration space is $2^{n-1}n!(n+1)$, which is the number of connected components.
\end{conjecture}

A careful analysis of which dihedral orders in $\M_{0,2(n-1)-3}$ contribute a connected component in $\M_{C_n}$ after folding, together with Theorem~\ref{thm:mainResult}, should provide a proof. For a tropical approach to the configuration space $\M_{C_n}$ see \cite{tropbin}.

\subsection{Commutative algebra of $\mathbf{u}$-relations} \label{subsec: commalg}
One can define a partial compactification of $\M_{0,n}$ as the closure of the dihedral embedding. Concretely, recall the $u$-relations from Section~\ref{subsec: 1.3}. These cut out $\M_{0,n}$ set theoretically as a very affine variety. Then, we may define $\widetilde{\M}_{0,n}:=\mathbb{V}(R_{ij})$ as the vanishing set of the $u$-relations in affine space $\CC^{\binom{n}{2}-n}$.
The semi-algebraic set $\widetilde{\M}_{0,n}(\RR_{>0})$ is a curvy version of the associahedron $A_{n-3}$ --
its boundary stratification is isomorphic to the face lattice of $A_{n-3}$~\cite{Brown09}. Defining $u$-relations from other simplicial complexes gives rise to other positive geometries with interesting stratifications; these have have been studied as \emph{binary geometries} in~\cite{pellytope, tropbin, delpezzo, HeLiRamanZhang}.

\smallskip

One motivation to consider this affine variety is its relation to the DMK compactification. The partial compactifications $\widetilde{\M}_{0,n}^\alpha$ for all dihedral charts is a set of open charts covering $\overline{\M}_{0,n}$. Indeed, each point on the boundary $\overline{\M}_{0,n} \setminus \M_{0,n}$ represents a \textit{nodal} stable curve $C$: a rational nodal curve with the $n$ marked points distributed across the irreducible components in such a way that each component contains no less than three special points, where special points are marked points or nodes. The tropicalization of $C$ is a tree with a vertex for each connected component, an edge for each node, and an unbounded edge for each of the $n$ marked points. Choosing an embedding of this tropical curve as a planar graph determines a cyclic ordering of the unbounded edges, which defines a dihedral ordering $\alpha$ of $[n]$. Then one can check that $C$ is contained in $\widetilde{\M}_{0,n}^\alpha$. As a limit of smooth curves in $\M_{0,n}$, a boundary curve represents the collision of two or more of the marked points - the colliding points ``bubble off" into a new component, with their configuration on this new component parametrizing their rate of approach with respect to each other.

\medskip

A natural commutative algebra question to ask is whether the primitive $u$-relations generate a prime ideal. 
This follows indirectly from~\cite[Theorem 3.3]{ClusterConfigSpaces}, which in particular says that $\widetilde{\M}_{0,n}$ is a smooth and irreducible affine scheme. 
Arkani-Hamed, He and Lam prove this by showing that $\widetilde{\M}_{0,n}$ is isomorphic to an affine open inside a smooth projective toric variety. 
The recursive structure of the boundary $\widetilde{\M}_{0,n} \setminus \M_{0,n}$ allows for an inductive argument on $n$.

We checked directly, up to $n=7$, that the ideal $I_n \subset \CC[u_{ij}]$ generated by the $u$-relations for $\M_{0,n}$ is prime.  
In \cite{M2Code}, the reader can find the code to generate $I_n$ and verify that \texttt{isPrime} terminates for $n=5,6,7$.
We propose a number of strategies to address the following problem. They will all boil down to a hard Gr\"{o}bner basis computation.

\begin{problem}
    Find a direct proof that $I_n$ is prime for all $n\geq 8$.
\end{problem}

An inductive argument would be to use the map $f:\widetilde{\mathcal{M}}_{0,n+1} \to \widetilde{\mathcal{M}}_{0,n}$ which forgets one of the marked points. This is a monomial map in the $u$-variables which is flat, and hence the proof reduces to checking $(f^*)^{-1}(I_{n+1}) = I_n$. A second strategy is to use the Jacobian criterion directly on the $u$-relations to prove $\mathbb{V}(I_n)$ is smooth as an affine scheme. A third approach is to directly compute a Gröbner basis of $\langle R_{ij}\rangle \subset \CC[u_{ij}^\pm]$ such that eliminating $u_{ij}^{-1}$ yields $I_n$. For $n=5$, the five $u$-relations are a Gröbner basis, but for $n=6$ the na\"{i}ve choices of elimination orders yield Gröbner bases that contain more than just the (extended) $u$-relations.
Finally, using~\cite[Proposition 23]{GSS}, we could iteratively eliminate variables until we reach an ideal for which proving primality is easier.
We have a candidate set $S$ of $u$-variables for which the elimination ideal of $I_n$ is a principal ideal generated by one of the primitive $u$-relations $R_{ij}$. Concretely, we conjecture that we can take $S$ to be the set of all $u$-variables that do not appear in $R_{ij}$, for $u_{ij}$ a longest chord in the $n$-gon.

\vfill

\noindent {\bf Acknowledgement}: 
We thank Hadleigh Frost for first showing us~\cite[Conjecture 11.1]{ClusterConfigSpaces} and Bernd Sturmfels for very useful feedback. We also thank Thomas Lam for pointing us to~\cite[Theorem 3.3]{ClusterConfigSpaces}, and Maximilian Wiesmann for suggesting~\cite[Proposition 23]{GSS}, both mentioned in Section~\ref{subsec: commalg}.

\bibliographystyle{plainurl}
\bibliography{references}

\end{document}